\documentclass{amsart}
\usepackage{amssymb}
\usepackage{amsfonts}
\usepackage{amsmath}
\usepackage[legalpaper,bookmarks=true,colorlinks=true,linkcolor=blue,citecolor=blue]{hyperref}
\usepackage{graphicx}%
\usepackage{fancyhdr}
\usepackage{color}
\usepackage[mathlines]{lineno}
\usepackage{lscape}
\usepackage{epsfig}
\usepackage{natbib}
\usepackage{geometry}
\usepackage{tgbonum}
\usepackage[]{listings}

\fontfamily{qcr}\selectfont

\newtheorem{theorem}{Theorem}
\theoremstyle{plain}

\newtheorem{corollary}{Corollary}

\newtheorem{proposition}{Proposition}

\numberwithin{equation}{section}

\begin{document}
\Large
\title[The exact probability law for the approximated similarity from the Minhashing method]{The exact probability law for the approximated similarity from the Minhashing method}

\author{Soumaila Dembele}
\author{Gane Samb Lo}

\begin{abstract} We propose a probabilistic setting in which we study the probability law of the Rajaraman and Ullman \textit{RU} algorithm and a modified version of it denoted by \textit{RUM}. These algorithms aim at estimating the similarity index between huge texts in the context of the web. We give a foundation of this method by showing, in the ideal case of carefully chosen probability laws, the exact similarity is the mathematical expectation of the random similarity provided by the algorithm. Some extensions are given.\\

\noindent \textbf{R\'{e}sum\'{e}.} Nous proposons un cadre probabilistique dans lequel nous \'{e}tudions la loi de probabilit\'{e} de l'algorithme de Rajaraman et Ullman \textit{RU} ainsi qu'une version modifi\'{e}e de cet algorithme not\'{e}e \textit{RUM}. Ces alogrithmes visent \`{a} estimer l'indice de la similarit\'{e} entre des textes de grandes tailles dans le contexte du Web. Nous donnons une base de validit\'e de cette m\'{e}thode en montrant que pour des lois de probabilit\'{e}s minutieusement choisies, la similarit\'{e} exacte est l'esp\'{e}rance math\'{e}matique de la similarit\'{e} al\'{e}atoire donn\'{e}e par l'algorithme \textit{RUM}. Des g\'en\'eralisations sont abord\'ees.\\

\noindent $^{\dag}$ Soumaila Dembele\\ 
Universit\'e des Sciences Sociale et de Gestion de Bamako ( USSGB)\\
Facult\'e des Sciences \'Economiques et de Gestion (FSEG)\\
Email: soumaila.demebele@ugb.edu.sn\\

\noindent $^{\dag}$ $^{\dag}$ Gane Samb Lo.\\
LERSTAD, Gaston Berger University, Saint-Louis, S\'en\'egal (main affiliation).\newline
LSTA, Pierre and Marie Curie University, Paris VI, France.\newline
AUST - African University of Sciences and Technology, Abuja, Nigeria\\
gane-samb.lo@edu.ugb.sn, gslo@aust.edu.ng, ganesamblo@ganesamblo.net\\
Permanent address : 1178 Evanston Dr NW T3P 0J9,Calgary, Alberta, Canada.\\

\end{abstract}
\maketitle

\newpage
\section{Introduction} \label{sec1}

\noindent In this paper, we are concerned with the evaluation of an important algorithm destined to provide the approximation of the exact similarity of two texts, in the frame of Web mining. \cite{AnandJeffrey} proposed a detailed algorithm we denote here as the \textit{RU} one. This algorithm is based on  \textit{minhashing} methods. To fix the ideas, let us consider two sets $S_{1}$ and $S_{2}$, whose total cardinality is $n$. The Jaccard similarity
between $S_{1}$ and $S_{2}$ is defined by:

\begin{equation}
p\mathbb{=}\frac{\#(S_{1}\cap S_{2})}{\#(S_{1}\cup S_{2})}  \label{pur1}
\end{equation}

\noindent Although this expression is simple, its computation is extremely time consuming in the context in Web mining, where the data may be huge.
For this reason, approximations based on probability theory and statistical methods are used.\newline

\noindent Before we come back the our precise subject, it may be useful to say some words on the general matter. The concept of similarity has been studied and is still studied by researchers from a variety of disciplines: (see e.g. \citealp{SteinEissen06}, \citealp{GionisIndykyMotwaniz}, for visual similarity, \citealp{Sung07} for the use of density functions in similarity detection, \citealp{GowerLegendre} and \citealp{ZezulaAmatoDohnal06} for the metric space approach, \citealp{StrehlGhoshMooney00}, \citealp{Formica05} in the context of information sciences, \citealp{BilenkoMooney03} and  \citealp{Theobald} for focus similarity on large-web collections).\medskip

\noindent The current work uses a \textit{minhashing} method (see e.g. \citealp{Fran}) on the Iterative Universal Hash Function Generator for Minhashing, and the
resemblance and containment of documents (see e.g. \citealp{Broder}). \newline

\noindent One way to deal with such problems is to transform the data into a low dimension representation, supposed to preserve enough information, and to derive the similarity index on the transformed data. In order to reduce the dimensionality of a data set, some  methods consist of introducing variables and feature selections or, of using a probabilistic dimension reduction technique (see e.g. \citealp{Guyon},  \citealp{Guyon1}, \citealp{Lawrence}, etc). The method we work on it in this papers uses the technique of \textit{signatures}. Let us explain this.\\

\noindent Consider $m$ subsets of $S$: $S_{1},...,S_{m}$. These sets can be represented as in Table \ref{tab21}, that we will call the \textit{representation matrix} or simply the \textit{signature}  of $S_{1},...,S_{m}$. This representation is set up as follows.

\begin{table}[h]
	\centering
	\begin{tabular}{|l|l|l|l|l|l|l|l|l|}
		\hline
		Elements & $S_{1}$ & $S_{2}$ & ... & $S_{h}$ & ... & $S_{l}$ & ... & $S_{m}$
		\\ \hline
		$1$ & $1$ & $0$ & $...$ & $0$ & $...$ & $1$ & $...$ & $1$ \\ \hline
		$2$ & $0$ & $0$ & $...$ & $1$ & $...$ & $0$ & $...$ & $0$ \\ \hline
		$...$ & $0$ & $...$ & $...$ & $...$ & $...$ & $...$ & $...$ & $...$ \\ \hline
		$i$ & $1$ & $0$ & $...$ & $1$ & $...$ & $1$ & $...$ & $1$ \\ \hline
		$...$ & $...$ & $...$ & $...$ & $...$ & $...$ & $...$ & $...$ & $...$ \\ 
		\hline
		$n$ & $0$ & $0$ & $...$ & $0$ & $...$ & $0$ & $...$ & $1$ \\ \hline
	\end{tabular}
	\caption{Representation matrix of $S_1$, \ldots, $S_m$}
	\label{tab21}
\end{table}

\begin{itemize}
	\item We form a rectangular array of $m+1$ columns.
	
	\item We put $S,S_{1},...,S_{m}$ in the first row.
	
	\item We put in the column of $S$ all the elements of $S$, that we might
	write from $1$ to $n$ in an arbitrary order.
	
	\item In each column $S_{h}$, $1\leq \ell \leq m$, we will put $1$ or $0$ on the row $i$
	depending on whether the $i^{th}$ element of $S$ is in $S_{h}$ or not.\\
\end{itemize}

\noindent It is immediate from Table \ref{tab21} that the following properties hold, for each couple $(i,j)$ such that $1\leq h \neq \ell\leq m$ :\\

\noindent (a) the cardinality of  $(S_{h}\cup S_{\ell})$ is the number of rows in Table \ref{tab21} crossing columns $S_h$ and $S_{\ell}$ at least with a unity value.

\noindent (b) the cardinality of $(S_{h}\cap S_{\ell})$ is the number of rows in Table \ref{tab21} crossing both columns $S_h$ and $S_{\ell}$ with a unity value.\\

\bigskip \noindent Hence, the representation matrix allow to get, visually, the similarity between $S_h$ and $S_{\ell}$. In particular, by denoting 
$S_{h}=(S_{ih}, {1\leq i\leq n)^{T}}$ for $1\leq h \leq m$, where $^{T}$ stands for the transpose of a matrix, we have

\begin{equation}
sim(S_{h},S_{\ell})=\frac{\#\{i,1\leq i\leq n,S_{ih}=S_{i \ell}=1\}}{\#\{i,1\leq i\leq n,(S_{ih}+S_{i \ell}=1)+(S_{ih}=S_{i \ell}=1)\}}, \label{pur1a}
\end{equation}

\noindent which can be written as

\begin{equation}
sim(S_{h},S_{\ell})=\frac{\#\{i,1\leq i\leq n,S_{ih}+S_{i \ell}=2\}}{\#\{i,1\leq i\leq n,(S_{ih}+S_{i \ell}=1)+(S_{ih}+S_{i \ell}=2)\}}. \label{pur1b}
\end{equation}

\noindent The \textit{RU} algorithm consists of reducing this matrix to a much less one with the help of \textit{minhashing} functions and, of computing the similarity between two columns which is meant to approximate the similarity between the sets represented by these columns.\\

\noindent Although it may be empirical observed that the approximation may be relevant, it does not exist, up to our knowledge, an theoretical evaluation of the discrepancy between the  exact similarity and the estimated similarity provided by the \textit{RU} method. The papers will fill this gap. Beyond that, it lays out a probabilistic frame to handle the problem and opens new research trends.\\

\noindent The rest of the paper is organized as follows. In Section \ref{sec2}, the full description of the \textit{RU} algorithm is given. Interesting remarks and properties will be addressed. Computation aspects will also be highlighted in this section. New forms, more appropriate to address the probability problem, will be given and a slightly modified algorithm, named \textit{RUM}, is proposed. In section \ref{sec3}, the \textit{RU} algorithm will be approached in a probability theory frame and the probability law of the random similarity is given and some consequences, among them the deviation from the true similarity, is characterized. Next, we study some conditions under which convergence of the \textit{RU} is explained. Finally, in a pure and random scheme, we completely justify the \textit{RU} algorithm in the probabilistic approach. Concluding remarks are stated in Section \ref{sec4}.

\section{RU and RUM algorithms} \label{sec2}

It is based on the notion of minhashing to reduce sets of huge sizes into sets of small sizes called signatures. The computation of the similarity is done on their compressed versions, i.e, on their signatures. To better explain this notion, let us consider $m$ subsets of a reference set $S$ of size $n\geq 1$ and let us use their \textit{representation matrix} as in Table \ref{tab21}.\\

\noindent Let us consider $k\geq 1$ functions $z_{\alpha}$ $(\alpha=1,...,k)$ from $\{1,...,n\}$ to itself in the following form:

\bigskip

\begin{equation}
z_{\alpha}(i)=x_{\alpha} i+y_{\alpha}\text{ mod }n,  \label{cong}
\end{equation}

\bigskip

\noindent where $x_{\alpha}$ and $y_{\alpha}$, $1\leq \alpha \leq k$, are given integers. We modify this function in the following way: $z_{\alpha}(i)=n$ when the remainder of the Euclidean division is zero.

\bigskip

\noindent Next, we extend Table \ref{tab21} by adding $k$ columns $\mathcal{Z}_{1},...,\mathcal{Z}_{k}$, such that each $\mathcal{Z}_{i}$ is the transpose of $(z_{\alpha}(1),...,z_{\alpha}(n))$. The resulting table is Table \ref{tab52}. Let us denote by $\mathcal{Z}$ the $(n \times k)$-matrix whose columns are $\mathcal{Z}_{1},\mathcal{Z}_{2},...,\mathcal{Z}_{k}$ and let us denote its lines rows by $\mathcal{Z}^{(1)},\mathcal{Z}^{(2)},...,\mathcal{Z}^{(n)}$.\\

\noindent The \textit{RU} algorithm replaces Table \ref{tab21} by a shorter one called \textit{minhashing} signature represented in Table \ref{tab53}

\bigskip

\begin{table}[tbp]
	\centering
	\begin{tabular}{|l|l|l|l|l|l|l|l|l|}
		\hline
		Element & $S_{1}$ & $S_{2}$ & ... & $S_{m}$ & $\mathcal{Z}_{1}$ & $...$ & $\mathcal{Z}_{k}$ & 
		\\ \hline
		1 & 1 & 0 & ... & 0 & $z_{1}(1)$ & ... & $z_{k}(1)$ & $\mathcal{Z}^{(1)}$ \\ \hline
		2 & 0 & 0 & ... & 1 & $z_{1}(2)$ & ... & $z_{k}(2)$ & $\mathcal{Z}^{(2)}$ \\ \hline
		. & 0 & . & ... & . & . & ... & . & . \\ \hline
		i & 1 & 0 & ... & 1 & $z_{1}(i)$ & ... & $z_{k}(i)$ & $\mathcal{Z}^{(i)}$ \\ \hline
		. & . & . & ... & . & . & . & . & . \\ \hline
		$n$ & 0 & 0 & ... & 0 & $z_{1}(n)$ & ... & $z_{k}($n$)$ & $\mathcal{Z}^{(n)}$ \\ 
		\hline
	\end{tabular}%
	\caption{Extension of the representation matrix by minhashing columns}
	\label{tab52}
\end{table}

\bigskip

\begin{table}[tbp]
	\centering
	\begin{tabular}{|l|l|l|l|l|}
		\hline
		minhashes & $t(S_{1})$ & $t(S_{2})$ & $...$ & $t(S_{m})$ \\ \hline
		$1$ & $c_{11}$ & $c_{12}$ & $...$ & $c_{1m}$ \\ \hline
		$2$ & $c_{21}$ & $c_{22}$ & $...$ & $c_{2m}$ \\ \hline
		... & $...$ & $.$ & $...$ & $.$ \\ \hline
		$k$ & $c_{k1}$ & $c_{k3}$ & $...$ & $c_{km}$ \\ \hline
	\end{tabular}%
	\caption{ Signature matrix}
	\label{tab53}
\end{table}

\bigskip

\noindent Table \ref{tab53} is obtained as follows, according to the method described in \cite{AnandJeffrey}, page 65.\\

\noindent \textbf{Algorithm of filling the columns $S_{j}$}.

\begin{itemize}
	\item[1.] Set all the $c_{\alpha j}$ equal to $\infty .$
	
	\item[2.] For each column $S_{j},$ proceed like this
	
	\begin{itemize}
		\item[2-a.] for each element $i,$ from $1$ to $n$, compute $%
		z_{1}(i),z_{2}(i),........,z_{k}(i)$.
		
		\item[2-b.] if $i$ is not in $S_{j},$ then do nothing and go to $i+1$
		
		\item[2-c.] if $i$ is in $S_{j},$ replace all the rows $(c_{\alpha j})_{1\leq
			r\leq k} $ by the minimum: min($c_{\alpha j},h_{\alpha}(i)).$
		
		\item[2-d.] go to $i+1$
	\end{itemize}
	
	\item[3.] Go to $j+1$
	
	\item[4.] End.
\end{itemize}

\bigskip

\noindent At the end of the procedure, each column will contain only integers between $1$ and $n$. The estimated similarity between two subsets $S_{h}$ and $S_{\ell}$, $1\leq h \leq \neq \ell\leq m$, based on this compressed table, is taken as the similarity of the columns $t(S)_h$ and $t(S)_\ell$ in the signature matrix which is

\begin{equation}
simRU(S_{h},S_{\ell})=\frac{\#(t(S_{h}) \cap t(S_{\ell}))}{k}.  \label{pur2}
\end{equation}

\noindent The sets $t(S_h)$ and $t(S_{\ell})$ are subsets of $\{1, 2, ...,n\}$ and stand for transformed sets of $S_h$ and $S_{\ell}$ through the \textit{minhashing} procedure.\\

\noindent From there, it is very important to give this remark. The algorithm is meant to gain time and to get an approximation of the similarity. It is based on the representation matrix. But, if we spent the required time to get it, there is nothing else to do, since the exact similarity is automatically read in virtue of Formulas \ref{pur1a} and \ref{pur1b}. We have to modify the \textit{RU} algorithm form a practical point of view.\\

\noindent The resulting modification, called \textit{RUM}, consists of the following. Suppose that we want to find and estimated similarity between $S_h$ and $S_{\ell}$. We proceed as follows.

\begin{itemize}
	\item[1.] Form one set $S_{h,\ell}$ by putting the elements of $S_h$ and then the
	elements of $S_{\ell}$ by putting twice elements of the intersection.\\
	
	\item[1.] Form the representation matrix with $N=n_{1}+n_{2}$ lines.
	
	\item[2.] Apply the RU algorithm to this collection by using Criterion (\textit{C}).
\end{itemize}

\bigskip

\noindent We do not seek to find the intersection. Elements of the intersection are counted twice here. The result is that we do not loose time in forming the representation matrix.\\

\noindent But, we will now have two approximations. First, we replace the representation of the \textit{RU} approach by that of the \textit{RUM} one. Next, we replace the latter by the signature matrix.\\

\noindent How is affected the original similarity? The \textit{RUM} algorithm actually seeks at estimating the modified similarity between two subsets $S_h$ and $S_\ell$, $1\leq h,\ell\leq m$ . Let us denoted by $simM(S_h, S_{\ell})$. It is immediately seen that we  still have a zero similarity index, that is $simM(S_h, S_{\ell})=0$, if the two sets $S_{h}$ and $S_{\ell}$ are disjoint, and a $100\%$ index if the sets are identical. In the general case, the number of rows is now $\#(S_h \cup S_{\ell})+\#(S_h \cap S_{\ell})$ and the common elements of the columns $S_h$ and $S_{\ell}$ is $\#(S_h \cap S_{\ell})$. The modified similarity between $S_h$ and $S_{\ell}$,  is 

\begin{equation}
simM(S_h \cap S_{\ell})=\frac{2\#(S_h \cap S_{\ell})}{\#(S_h \cup S_{\ell})+\#(S_h \cap S_{\ell})},
\end{equation}

\noindent which gives

\begin{equation} \label{simModD}
simM(S_h \cap S_{\ell})=\frac{2 sim(S_h \cap S_{\ell})}{1+sim(S_h \cap S_{\ell})}
\end{equation}

\noindent and, reversely, 

\begin{equation} \label{simModI}
sim(S_h \cap S_{\ell})=\frac{simM(S_h \cap S_{\ell})}{2-simM(S_h \cap S_{\ell})}
\end{equation}

\noindent It is also clear that from the previous formulas that $sim$ and $simM$ take any of the value zero and one simultaneously.\\

\noindent We adopt the following rule : We use the modified \textit{RU} algorithm in place of the original one. We will avoid to find the intersection, which in fact would stop our procedure since the similarity is already found, and by then, we gain a huge amount of time. At the end of the \textit{RU} algorithm implementation on the modified set, we apply Formula \ref{simModI}, to get the approximation

\begin{equation} \label{simModRUM}
simRU(S_h \cap S_{\ell})=\frac{simRUM(S_h \cap S_{\ell})}{2-simRUM(S_h \cap S_{\ell})}
\end{equation}

\noindent Now, the question is how accurate is the approximation? Empirical studies strongly support the method. For instance, the four canonical Gospels have been compared with the target of assessing the hypothesis of the existence of a hidden or lost sources, named \textbf{Q} source, from which the current gospels are derived. The Gospel have been transformed into sets of words of $p=3$ letters (named $p$-\textit{shingles}). The numbers of $3$-shingles of the four gospels are at least 55.000 and at most 110.000. The shortest time to compute the exact similarity between two gospels is around eight (8) minutes while the computation of the similarity between John and Matthews Gospels requires 3080 seconds (around 51 minutes). By using the \textit{RUM} algorithm with only a small number $k=5$ of \textit{minhashing} functions, estimations of the similarity indices are obtained with a much smaller time, around 20 seconds. The estimated values showed clear trends for the exact values.\\

\noindent Clearly here, the choice of the \textit{minhashing} functions is arbitrary. Without saying it, their choice is subject to a probability law. Implicitly, the uniform law is assumed. Even in that implicit choice, we did not see a study on the exact final probability law.\\

\noindent In the forthcoming section, we deal with the probability law of \textit{simRU}, considered as a random variable. We will study it with respect to the probability law of the random coefficients 
$a_i$ and $b_i$, $i=1,...,k$.\\

\section{Probabilistic approach} \label{sec3}

\noindent Let us give a probabilistic approach of the similarity.\\

\noindent \textbf{I - The similarity as a conditional probability}.\\

\noindent We adopt the notation introduced in the previous section, in particular the representation matrix in Section \ref{tab21}. Now, we suppose that we pick at random 
a row $X$ from the number of lines in Column one in the Table \ref{tab21} and for each $h$, $1\leq h\leq m$, let $X_{k,\ell}$ be the Bernoulli random variable taking the value at the crossing between the column $S_h$ and the row $X$. We are going to see that the random variable $X$ guide the similarity index.

\begin{theorem} \label{t01} Let us randomly pick a row $X$ among $n$ rows. Let $S_{X,h}$
	be the value of the row $X$ at the crossing with a column $S_h$, $1\leq h\leq m$ in the representation matrix. Then the
	similarity between two sets $S_{h}$and $S_{\ell}$, $1\leq h \leq \ell$,  is the
	conditional probability of the event $(S_{X,h}=S_{X,\ell}=1)$ with respect to
	the event $(S_{X,h} + S_{X,k}\geq 1)$. i.e 
	\begin{equation*}
	sim(S_{\ell},S_{h})=\mathbb{P}[(S_{X,h}=S_{X,\ell}=1)/(S_{X,h}+S_{X,\ell}\geq
	1)].
	\end{equation*}
\end{theorem}

\noindent \textbf{Proof.} We first observe that for the defined matrix
below, the set of rows can be split into three classes, based on the columns 
$S_{\ell}$ and $S_{h}$:\newline

\noindent 1. The rows $(A)$ that cross both columns $S_{\ell}$ and $S_{h}$ with unity values.\newline

\noindent 2. The rows $(B)$ that cross $S_{\ell}$ and $S_{h}$ with a unity value and a null value.\newline

\noindent 3. The rows $(C)$ that cross both $S_{\ell}$ and $S_{h}$ with null values.\newline

\noindent Let us show that $sim(S_{\ell},S_{h})=\mathbb{P}[(S_{X,h}=S_{X,\ell}=1)/(S_{X,h} + S_{X,\ell}\geq 1)]$. \newline

\noindent Clearly, the similarity is the ratio of the number of rows $(A)$ to the sum of the numbers of rows $X$ and the number of rows $(B)$. The rows $(C)$
are not involved in the similarity between $S_{h}$ and $S_{k}$. Thus 

\begin{equation*}
sim(S_{\ell},S_{h})=\frac{\#\{i,1\leq i\leq n,S_{X,h}=1,S_{X,\ell}=1\}}{%
	\#\{i,1\leq i\leq n,(S_{X,h}+S_{X,\ell}=1)+(S_{X,h}=1,S_{X,\ell}=1)\}}.
\end{equation*}

\bigskip

\noindent Then, by dividing the numerator and the denominator by $n$, we
will have

\begin{equation*}
sim(S_{\ell},S_{h})=\frac{\frac{\#\{i,1\leq i\leq n,S_{X,h}=1, \ S_{X,\ell}=1\}%
	}{n}}{\frac{\#\{i,1\leq i\leq n,(S_{X,h}+S_{X,\ell}=1)+(S_{X,h}=1,S_{X,\ell}=1)\}}{n}}.
\end{equation*}

\noindent Hence we get the result

\begin{equation*}
sim(S_{\ell},S_{h})=\mathbb{P}[(S_{X,h}=S_{X,\ell}=1)/(S_{X,h} + S_{X,\ell}\geq 1)].
\end{equation*}

\bigskip

\noindent This theorem will be the foundation of the statistical estimation
of the similarity as a probability.\newline

\noindent \textbf{Important remark.} When we consider the similarity of two subsets, say $S_h$ and $S_k$ and we use the global space as $S_h \cup S_{\ell}$,
we may see that the similarity is, indeed, a probability. But when we simultaneously study the joint similarities of several subsets, say at least $S_h$, $S_{k}$ and $S_\ell$ with the global set $S_h \cup S_k \cup S_{\ell}$, the similarity between two subsets is a \textbf{conditional} probability. Then, using the fact that the similarity is a probability to prove the triangle inequality is not justified, as claimed in \cite{AnandJeffrey}, page 76.\\

\noindent \textbf{II - Expected or Normal Similarity}.\\

\noindent Before we begin, we stress that the notation $k$ and $m$ are not related to those in the other sections. These notation should stay specific to the problem handled here.\\

\noindent We shall use the language of the urns. Suppose that we have a reference set of size $n$ that we take as an urn \textbf{U}. We pick at random a subset $X$ of size $k$ and a subset $Y$ of size $m$. If $m$ and $k$ have not the same value, the picking order of the sets does have an impact on our results. We then proceed at the beginning by picking at random the
first subset, that will be picked all at once, next put it back in the urn \textbf{U} (reference set). Then we pick the other subset. Let us ask
ourselves the question : what is the expected value of the similarity of Jaccard?\newline

\noindent The answer at this question allows us later to appreciate the degree of similarity between the texts. We have the following result :%
\newline

\begin{proposition}
	Let $U$ be a set of size $n$. Let us randomly pick two subsets $X$ and $Y$
	of $U$, of respective sizes $m$ and $k$ according to the scheme described
	above. We have \begin{equation}
	\mathbb{P}(Card(X\cap Y)=j) 
	= \frac{1}{2}\left( \frac{C_{k}^{j} \ C_{n-k}^{m-j}}{C_{m}^{n}} +\frac{C_{m}^{j} \ C_{n-m}^{k-j}}{C_{k}^{n}}   \right) \mathbb{I}_{(0\leq j \leq \min(k,m))}.\label{sim01} 
	\end{equation}
	
	\noindent For all $p\leq 1$, the $p$-th moment of the random similarity $sim(X,Y)$ is given by 
	\begin{equation}
	\mathbb{E}(sim(X,Y)) 
	=\sum_{j=0}^{\min (k,m)}\frac{j}{2(m+k-j)}
	\left\{ \frac{C_{k}^{j} \ C_{n-k}^{m-j}}{C_{m}^{n}} +\frac{C_{m}^{j} \ C_{n-m}^{k-j}}{C_{k}^{n}}\right\}. \label{sim02}
	\end{equation}
\end{proposition}

\noindent \textbf{Proof}. Let us use the scheme described above. Let us
first pick the set $X$. We have $L=C_{n}^{k}$ possibilities. Let us denote
the subsets that would take $X$ by $X_{1},...,X_{L}.$ The searched
probability becomes
\begin{eqnarray*}
	\mathbb{P}(Card(X\cap Y)=j)&=&\sum_{s=1}^{L}\mathbb{P}((Card(X\cap Y)=j)\cap X_{s})\\
	&=&\sum_{s=1}^{L}\mathbb{P}((Card(X\cap Y)=j)/X_{s})\mathbb{P}(X_{s}).
\end{eqnarray*}

\noindent Once $X_{s}$ is chosen and fixed, we get
\begin{equation*}
\mathbb{P}((Card(X\cap Y)=j)/X_{s})=\frac{C_{m}^{j}\text{ \ \ }C_{n-m}^{k-j}%
}{C_{m}^{k}}\text{ }.
\end{equation*}

\noindent To explain this, we start with the fact that $X_{s}$ is fixed and contains $C^{j}_{k}$ combinations of $j$ elements. Now, we have to choose a combination $\mathcal{C}$ 
of $m$ elements from $n$ elements which contains one of the $C^{j}_{k}$ combinations of $X_{s}$ in such a way that none of the other elements of $\mathcal{C}$ is in $X_{s}$. This means that
on should choose first a combination of $j$ among the $k$ elements of $X_s$ with $C^{j}_{k}$ ways, and next one completes with a combination of $m-j$ elements among the $n-k$ elements of the complement of  $X_{s}$.\\

\noindent Now, since $\mathbb{P}(X_{s})=1/C_{n}^{k}=1/L$, we conclude  that

\begin{equation*}
\mathbb{P}(Card(X\cap Y)=j)=\sum_{s=1}^{L}\frac{C_{k}^{j} \ C_{n-k}^{m-j}}{C_{m}^{n}}(1/L)=\frac{C_{k}^{j} \ C_{n-k}^{m-j}}{C_{m}^{n}}.
\end{equation*}

\noindent The result corresponding to picking up $Y$ first, is obtained by symmetry of
roles of $k$ and $n$. We then get (\ref{sim01}). The formula (\ref{sim02})
comes out immediately since
\begin{equation}
sim(X,\text{\ }Y)=\frac{\#(X\cap Y)}{\#(X\cup Y)}=\frac{\#(X\cap Y)}{%
	m+k-\#(X\cap Y)}.  \label{esp}
\end{equation}

\noindent $\blacksquare$\\

\bigskip \noindent \textbf{III - Approximated Similarity Based on the Strong Law of Large Number}.\\

\noindent Since the similarity is a conditional probability in according to Theorem $\ref{t01}$, we can deduce a strong law of Large numbers, which is by the way a Glivenko-Cantelli property in the discrete case, in the following way.

\begin{theorem} \label{gclaw}
	Let $sim(S_1,S_2)$ be the similarity between two subsets $S_1$ and $S_2$ of a set whose size is considered very large. Let us pick at
	random a subset $S_{1,n}$ from $S_1$ with size $n(1)$ and a subset $S_{2,n}$ from $S_2$ of size $n(2)$ and let us consider the random similarity $sim_{n}(S_1,S_2)=sim(S_{1,n},S_{2,n})$ between 
	$S_{1,n}$ and $S_{2,n}$. Then $sim_{n}(S_1,S_2)$ converges almost-surely to $sim(S_1,S_2)$ with at rate of convergence in the order of $(n(1)+n(2))^{-1/4}$ when $n(1)$ and $n(2)$ become
	simultaneously large.
\end{theorem}

\bigskip \noindent That is a direct consequence of the classical theorem of Glivenko-Cantelli.\\

\noindent Finally, we come to the probability law induced by the \textit{RU} algorithm.\\

\bigskip \noindent \textbf{IV - Probability law induced by the \textit{minhashing} method and application}.\\

\noindent \textbf{A - Two other alternative versions A simple criterion}.\\

\noindent Before to give two alternate versions of the \textit{RU/RUM} algorithm. The first will be particularly useful while addressing the probability law. Also it leads to a new procedure that will be the base of the implementation of the algorithm in computer packages.\\

\noindent \textbf{(a) A simple criterion}.\\

\noindent If we look carefully at the algorithm, we may see that we have the following criteria.\\

\noindent \textbf{Criterion (C)}. The transpose of each column 
$$
t(S_{h})=[(c_{\alpha h})_{1\leq \alpha \leq k}], \ 1\leq h \leq m,
$$

\noindent in Table \ref{tab53} is the minimum of the rows $\mathcal{Z}^{(i)}=(z_{1}(i),...,z_{k}(i))$ of Table \ref{tab52}, when $i$ covers the elements $i$ of $S_{h}$, where
the minimum is operated coordinate-wisely.\\

\bigskip \noindent The proof comes easily by looking at simple cases with small cardinalities. The induction to arbitrary cardinalities is immediate.\\

\noindent This simple remark allows to set up programs in a much easier way through a kind of Markov process.\\
\newline

\noindent \textbf{(b) A version if form a Markov process}.\\

\noindent We want to form the final transformed matrix signature as defined in Table \ref{tab53} by denoting $U_{h}$ as the column associated with $t(S_{h})$ and $U_{\ell}$ as the column associated with $t(S_{\ell})$. We remind that $U_{h}$ and $U_{\ell}$ are vectors of dimension $k$. This procedure will be implemented is easy to implement into computer packages.\\

\noindent We remind that in the original matrix signature, the elements of $S_{h}\cup S_{\ell}$ are given in an arbitrary order $(\sigma _{0}(i),$ $1\leq
i\leq N)$. We denote by $C(i,h)$ the value at which the row $i$ and the column $h$ cross each other in Table \ref{tab52}. In what follows, for any Boolean variable $C$, $\mathbb{I}(C)$ stands for
the indicator function of $C$, which takes the value one if $C$ holds and zero otherwise.\\

\noindent Let us express $RU$ algorithm as a final step of Markov process.\\

\noindent We fix $h$, $1\leq h\leq m$. The following procedure iteratively forms the final value of $U_h$.\\

\noindent Step $1$.  Do for each $h$, $1\leq h \leq m$ :
\noindent  $(1a)$.  Take $U_{h}^{0}=(n+1,n+1,...,n+1)^{t}\in \mathbb{R}^{k} $. Put $U_{h}^{0}$ in the column $t(S_{h})$ of Table \ref{tab53}.\\

\noindent Sub-step $(1b)$. For each $i$ from $1$ to $n$, we take
$$
U_{h}^{i}=U_{h}^{i-1} \mathbb{I}(C(i,h)=0)+\min ( U_{h}^{i-1},(\mathcal{Z}^{(i)})^{T}) \mathbb{I}(C(i,h)=1).
$$.

\noindent Step $2$.  For each $1\leq h \neq \ell \leq m$, compute the estimated similarity :
\begin{equation}
simRU(S_{h},S_{\ell},\sigma _{0})=\frac{1}{k}\sum_{\alpha=1}^{k}I(U_{h}(\alpha)=U_{\ell}(\alpha)). \label{estimSim}
\end{equation}

\noindent This second algorithm is more simple to implement.\\

\bigskip  \noindent \textbf{B - Probability laws}.\\

\noindent Now we are going to compute the estimated similarity \textit{simrum} in a complete randomly experience. We consider two subsets $S_{\ell}$ and $S_{\ell^{\prime}}$ of $S$. We allow
the elements of $S$ be ordered according to a permutation $\sigma$ of $\{1,2,...,n\}$. We denote set of permutations of $\{1,2,...,n\}$ as $\mathcal{S}_n$ and consider the probability space
$(\mathcal{S}_n, \mathcal{P}(\mathcal{S}_n), \mathbb{P}_{0})$, where $\mathcal{P}(\mathcal{S}_n)$ is the power set of $\mathcal{S}_n$ and $\mathbb{P}_{0}$ is the uniform probability measure on 
$\mathcal{S}_n$ defined by $\mathbb{P}_{0}(\{\sigma\})=\frac{1}{N!}$ for $\sigma \in \mathcal{S}_n$.\\

\noindent Next, we choose the following \textit{minhashing} function%
\begin{equation*}
Z_{\alpha}(i)=X_{\alpha} i+Y_{\alpha}\text{ mod }n, \ \ i=1,...,n.
\end{equation*}

\noindent with a random generation of the integers $(X_{1},Y_{1})$,\ldots, $(X_{k},Y_{k})$.\\

\noindent From there, a number of possibilities may be conceived. Do we take the $(X_{\alpha},Y_{\alpha})$'s as independent? independent and identically independent? dependent according a what copula? etc. We may also discuss about the dependence between $X_{\alpha}$ and $Y_{\alpha}$ for each $\alpha=1,\ldots,k$?\\

\noindent As a first step, let us suppose that :\\

\noindent (H) $(X_{1},Y_{1})$,\ldots, $(X_{k},Y_{k})$ are independent and identically distributed with common probability law $\mathbb{P}_{(X,Y)}$.\\

\noindent We apply the \textit{RUM} algorithm and observe the estimated random similarity between $S_{\ell}$ and $S_{\ell^{\prime}}$

\begin{equation*}
simrum(S_{\ell},S_{\ell^{\prime}})=\frac{1}{k}\sum_{\alpha=1}^{k}I(U_{\ell}(\alpha)=U_{\ell^{\prime}}(\alpha)).
\end{equation*}

\noindent If there is no risk of confusion, we simply write $simrum$ at the place $simrum(S_{\ell},S_{\ell^{\prime}})$ as we also use $sim$ and $simru$ at the place of $sim(S_{\ell},S_{\ell^{\prime}})$ and $simru(S_{\ell},S_{\ell^{\prime}})$ respectively. We are now going to give the probability law of $simrum$ after the following notations. The matrix $\mathcal{Z}$ in Table \ref{tab52} is random now and we denote

\begin{equation*}
\Gamma _{\ell }(\sigma )=\left\{ i\in \lbrack 1,n],\text{ \ }C(i,\ell
,\sigma )=1\right\}.
\end{equation*}

\noindent Introduce the following notation. Let $1\leq t\leq k$ define $\mathcal{B}_{t}$ as the set of all $t$-tuples. For $(\beta{1},...,\beta_{t})\in \mathcal{B}_{t}$, define 
$(\beta_{1},...,\beta_{t})^{c}$ as the complement of the set of $\left\{\beta_{1},...,\beta_{t}\right\} $ in $\{1,...,n\}$.\\

\noindent Define also for $1\leq p, q \leq n$, $1\leq \ell, \ell^{\prime} \leq m$, $\sigma \in \mathcal{S}$,

$$
\overline{m}(p,q,\sigma,\ell)=\min \{p+iq \mod{ } n, \ i\in \Gamma(\sigma) \}
$$

\noindent and

$$
B(\ell, \ell^{\prime}, \sigma)=\{ (p,q), 1\leq p, q\leq n,  \overline{m}(p,q,\sigma,\ell)=\overline{m}(p,q,\sigma,\ell^{\prime}) \}
$$

\noindent The probability law of the estimated similarity is the following.

\begin{theorem} \label{theo1} By Criterion (C) $U_{h}$ is given as follows.
	\begin{equation}
	U_{h}=\min (\mathcal{Z}^{(i)}, \ \ i\in \Gamma _{h}(\sigma ))^{T}\in \mathbb{R}^{k}, \ h =\ell, \ell^{\prime},  \label{eq1}
	\end{equation}
	
	\noindent where the minimum of rows is done by coordinate-wisely. Then $1\leq \alpha \leq k$,the probability of the event $(U_{\ell }(\alpha )=U_{\ell ^{\prime}}(\alpha ))$ is
	given by
	\begin{equation}
	\mathbb{P}(U_{\ell }(\alpha )=U_{\ell ^{\prime }}(\alpha ))=\frac{1}{n!}%
	\sum_{\sigma \in S_{n}}\sum_{(p,q)\in B(\ell ,\ell ^{\prime },\sigma )}%
	\mathbb{P}(X_{\alpha }=p,Y_{\alpha }=q).  \label{eq2}
	\end{equation}
	
	\noindent Moreover, the probability law of $simrum$ is given by the discrete probability measure defined on $\mathcal{V}=\{s\in [0,1], ks\in \mathbb{N}\}$ by  
	\begin{equation}
	\mathbb{P}(simrum=s)=\sum_{C\in B_{t}}\prod_{\alpha \in C}\mathbb{P}(\mathbb{(%
	}U_{\ell }(\alpha )=U_{\ell ^{\prime }}(\alpha ))\times \prod_{\alpha \notin
		C}\mathbb{P}(\mathbb{(}U_{\ell }(\alpha )\neq U_{\ell ^{\prime }}(\alpha )).
	\label{eq3}
	\end{equation}
	
	\noindent for $s\in \mathcal{V}$.
	
\end{theorem}

\noindent \textbf{Proof}. We are going to compute the probability of the event $(U_{\ell}(\alpha)=U_{\ell^{\prime}}(\alpha))$. First, it follows from the algorithm that (\ref{eq1}) is
straightforward, that is%
\begin{equation*}
U_{\ell }=\min (\mathcal{Z}^{(i)}, \ i\in \Gamma _{\ell}(\sigma))^{T}\in \mathbb{R}%
^{k}, \  \ell =1,2.
\end{equation*}

\noindent Now we are going to estimate the probability law of the event $(U_{\ell}(\alpha )=U_{\ell ^{\prime }}(\alpha))$. We point out that, conditionally on $\sigma $, the couples $(U_{\ell }(\alpha ),U_{\ell ^{\prime }}(\alpha ))$ are independent, and probabilities for events only depending on ($U_{\ell }(\alpha ),U_{\ell^{\prime }}(\alpha ))$, are computed with $\mathbb{P}_{(X_{\alpha
	},Y_{\alpha })}$. We get

\begin{eqnarray*}
	\hspace{1cm}&&\mathbb{P}(U_{\ell }(\alpha ) =U_{\ell ^{\prime }}(\alpha))\\
	&=&\mathbb{P}(\min (Z_{\alpha }(i),\ \ i\in \Gamma _{\ell }(\sigma ))=\min (Z_{\alpha }(i),\ \ i\in \Gamma _{\ell ^{\prime }}(\sigma )) \\
	&=&\mathbb{P(}\min (X_{\alpha}+iY_{\alpha }\text{ mod } N,\ \ i\in \Gamma_{\ell }(\sigma ))=\min (X_{\alpha }+iY_{\alpha }\text{ mod } N,\ \ i\in \Gamma _{\ell ^{\prime }}(\sigma )).
\end{eqnarray*}

\noindent By using the notation $\overline{m}(p,q,\sigma ,\ell)$ and $B(\ell ,\ell^{\prime},\sigma)$ introduced above, we finally get

\begin{eqnarray*}
	\hspace{3cm}\mathbb{P}(U_{\ell }(\alpha ) &=&U_{\ell ^{\prime }}(\alpha ))=\mathbb{P}_{(\sigma ,X_{\alpha },Y_{\alpha })}(B(\ell ,\ell ^{\prime },\sigma )) \\
	&=&\mathbb{P}_{\sigma }\otimes \mathbb{P}_{(X_{\alpha },Y_{\alpha })}(B(\ell
	,\ell ^{\prime },\sigma )) \\
	&=&\frac{1}{n!}\sum_{\sigma \in S_{n}}\mathbb{P}_{(X_{\alpha },Y_{\alpha
		})}(B(\ell ,\ell ^{\prime },\sigma )) \\
	&=&\frac{1}{n!}\sum_{\sigma \in S_{n}}\sum_{(p,q)\in B(\ell ,\ell ^{\prime
		},\sigma )}\mathbb{P(}X_{\alpha }=p,Y_{\alpha }=q)
\end{eqnarray*}

\noindent and (\ref{eq2}) is proved. Next, we recall that 
\begin{equation*}
simrum=\frac{\#\left\{ 1\leq \alpha \leq k,\text{ \ }U_{\ell }(\alpha
	)=U_{\ell ^{\prime }}(\alpha )\right\} }{k},
\end{equation*}

\noindent which entails
\begin{equation*}
simrum \in \#\left\{ \frac{t}{k},\text{ }1\leq \alpha \leq k\right\} ,\text{ }%
s=\frac{t}{k},\text{ }t=sk\in \mathbb{N}.
\end{equation*}

\bigskip \noindent Hence,

\begin{eqnarray*}
	\hspace{2cm}\mathbb{P}(simrum &=&s)=\mathbb{P}(simrum=t) \\
	&=&\mathbb{P}\left( \frac{\#\left\{ 1\leq \alpha \leq k,\text{ \ }U_{\ell }(\alpha
		)=U_{\ell ^{\prime }}(\alpha )=t\right\} }{k}\right).
\end{eqnarray*}

\noindent Let us put
\begin{equation*}
B_{t}=\{\alpha _{i_{1}}<\alpha _{i_{2}}<...<\alpha _{i_{t}}:U_{\ell }(\alpha
_{i})=U_{\ell ^{\prime }}(\alpha _{i}),\text{ \ \ }\forall 1\leq i\leq
t,U_{\ell }(r)\neq U_{\ell ^{\prime }}(r),\text{ }\forall r\in
\{1,...,k\},r\neq \alpha _{i}\}.
\end{equation*}

\noindent Let $D_{t}$ be the ordered subsets of $\{1,...,k\}$ of size $t$. If $(\alpha
_{i_{1}},\alpha _{i_{2}},...,\alpha _{i_{t}})\in D_{t}$, we denote  $(\alpha
_{i_{1}},\alpha _{i_{2}},...,\alpha _{i_{t}})^{c}$ as its complement in $%
D_{t}$.\\

\noindent Thus, we have
\begin{equation*}
B_{t}=\{c=(\alpha _{i_{1}},\alpha _{i_{2}},...,\alpha _{i_{t}})\in
D_{t}:\forall 1\leq i\leq t,\text{ }U_{\ell }(\alpha _{i})=U_{\ell ^{\prime
}}(\alpha _{i}),\text{ }\forall r\in c^{c}\text{\ \ },U_{\ell }(r)\neq
U_{\ell ^{\prime }}(r)\}.
\end{equation*}

\noindent Let $B\in B_{t}$, we have

\begin{eqnarray*}
	\hspace{2cm}B(c) &=&\bigcap_{\alpha \in c}(U_{\ell }(\alpha )=U_{\ell ^{\prime }}(\alpha
	))\bigcap \bigcap_{\alpha \notin c}(U_{\ell }(\alpha )\neq U_{\ell ^{\prime
	}}(\alpha )) \\
	&=&B_{1}(c)\bigcap B_{2}(c).
\end{eqnarray*}

\noindent We conclude that

\begin{eqnarray*}
	\hspace{2cm}\mathbb{P}(simrum &=&s)=\sum_{c\in B_{t}}\mathbb{P(}B_{1}(c)\bigcap B_{2}(c)) \\
	&=&\sum_{c\in B_{t}}\prod_{\alpha \in c}\mathbb{P(}U_{\ell }(\alpha )=U_{\ell ^{\prime }}(\alpha ))\times \prod_{\alpha \notin c}\mathbb{P}(U_{\ell }(\alpha )\neq U_{\ell ^{\prime }}(\alpha ))
\end{eqnarray*}

\noindent This completes the proof of Theorem.$\blacksquare$\\

\bigskip \noindent We have the following consequence.

\begin{corollary}\label{cor1} For any $p\in[0,1]$, we have 
	
	\begin{eqnarray}
	\mathbb{P}(|simrum-p| \leq\varepsilon )&=&\sum_{p-\varepsilon \leq s\leq p+\varepsilon }\mathbb{P}(simrum=s) \label{SimDev}\\
	&=&\sum_{[k(p-\varepsilon )\leq t\leq k(p+\varepsilon )]}\mathbb{P}(simrum=\frac{t}{k}). \notag
	\end{eqnarray}
\end{corollary}

\noindent \textbf{Proof}. Let us put $p=s$ in (\ref{pur1}).\newline

\noindent We get

\begin{eqnarray*}
	\mathbb{P}(|simrum-p| \leq\varepsilon )&=&\mathbb{P}(p-\varepsilon \leq
	simrum\leq p+\varepsilon ) \\
	&=&\sum_{p-\varepsilon \leq s\leq p+\varepsilon }\mathbb{P}(p_{3}=s).
\end{eqnarray*}

\noindent Since $s=t/k$, we obtain

\begin{equation*}
\mathbb{P}(|simrum-p|<\varepsilon )=\sum_{[k(p-\varepsilon )\leq t\leq
	k(p+\varepsilon )]}\mathbb{P}(simrum=\frac{t}{k}).
\end{equation*}

$\square$\\

\bigskip \noindent As a first application, Formula \ref{eq3} allows to compute the $p$-th moment of $simrum$, for $p\geq 1$, which is

$$fi
\mathbb{E}(simrum^p)=\sum_{s\in \mathcal{V}} s^p \mathbb{P}(simrum=s).
$$ 

\noindent Let us denote by SIMRUM the mathematical expectation of $simrum$, that is $SIMRUM=\mathbb{E}(simrum)$. In the context of discrete random variables, we may use Formula \ref{SimDev} to find a $95\%$-confidence interval by using an iterative procedure.\\

\noindent In conclusion, this result allows to have relevant confidence intervals of random modified similarity, and by then, of the random similarity through Formula \ref{simModRUM}. The comparison between the true similarity and the estimated similarity will no longer be done with a sole observation, but with respect to the whole confidence interval. This makes the comparison more reliable.\\

\noindent Several interesting questions remain open. For example, what is the efficiency of the estimation? What is the impact of the probability law of $\mathbb{R}^{2k}$-random variable

$$
((X_1,Y_1), \ldots, (X_k,Y_k)),
$$   

\noindent on the quality of the estimation? What happens as $k$ gets bigger?\\

\noindent Answering all these questions are beyond the scope of this paper. But, at least, we are going to give definite results on the \textit{RU} algorithm as a statistical method and lay out the general case.\\

\bigskip \noindent \textbf{C - Assessment of the \textbf{RU} algorithm as an estimation method}.\\

\noindent Let us begin to give the main idea of the method. This time suppose that the set $S=\{1, ...,n\}$ is given in a fixed order in the representation matrix and we write it from yhe top to the bottom in its natural order. It is attempted to consider $k$ random and independent permutations $Z_{\alpha}$, $1 \leq \alpha \leq k$ of the set $S=\{1, ...,n\}$. Suppose we are able to get them. We may complete the algorithm. At the arrival, we have the probability law of $simrum$ through the following notation. Let us introduce this new notation, for $L\subset \{1, 2, ...,n\}$, $\sigma \in \mathcal{S}_n$,

$$
Min_{L}(\sigma)=\min_{i\in L} \sigma(i),
$$

\noindent For $1\leq \ell \leq n$, we make the following abuse of notation and write $Min_{\Gamma(\ell)}=Min_{\ell}$. Now, denote for $1\leq \ell \leq n$, $1\leq \ell \leq \ell^{\prime}\leq m$,

$$
B(\ell, \ell^{\prime})=\{ \sigma \in \mathcal{S}_n, \pi_{\ell}(\sigma)=\pi_{\ell^{\prime}}(\sigma) \}.
$$

\noindent The probability law $simrum$ is still given by

\begin{equation*}
\mathbb{P}(simrum =s)\sum_{c\in B_{t}}\prod_{\alpha \in c}\mathbb{P(}U_{\ell }(\alpha )=U_{\ell ^{\prime }}(\alpha ))\times \prod_{\alpha \notin c}\mathbb{P}(U_{\ell }(\alpha )\neq U_{\ell ^{\prime }}(\alpha ))
\end{equation*}

\noindent with  

$$
\mathbb{P}(U_{\ell }(\alpha)=U_{\ell ^{\prime }}(\alpha ))=\mathbb{P}_{Z_{\alpha}}(B(\ell, \ell^{\prime})).
$$

\noindent for $ks \in \mathbb{N}$.\\

\noindent Here, what is expected is that two different rows $\mathcal{Z}^{i}$ and $\mathcal{Z}^{j}$, $1\leq i \leq j \leq$ will be probably disjoint, or at the least, that the probability they are not disjoints is very low. Suppose for a while that this is the case. Let us denote by $I(\ell-\ell^{\prime})=:I_{1}$ the set of lines pertaining to elements of $S_{\ell} \setminus S_{\ell^{\prime}}$,
$I(\ell^{\prime}-\ell)=:I_{3}$ the set of lines pertaining to elements of $S_{\ell^{\prime}} \setminus S_{\ell}$ and $I(\ell+\ell^{\prime})=:I_{2}$ the set of lines pertaining to elements of 
$S_{\ell} \cap S_{\ell^{\prime}}$. It is clear that

$$
U_{\ell}=\min \left(\min_{i \in I_2} \mathcal{Z}^{i}, \ \min_{i \in I_1} \mathcal{Z}^{i}\right)
$$

\noindent and

$$
U_{\ell^{\prime}}=\min \left( \min_{i \in I_2} \mathcal{Z}^{i}, \ \min_{i \in I_3} \mathcal{Z}^{i}\right).
$$

\noindent Thus, in the hypothesis of disjoint lines $\mathcal{Z}^{i}$'s, any event $BU(\alpha)=(U_{\ell}(\alpha)=U_{\ell^{\prime}}(\alpha))$ is surely achieved through the part

$$
\min_{i \in I_2} \mathcal{Z}^{i},
$$

\noindent meaning that $(U_{\ell}(\alpha)=U_{\ell^{\prime}}(\alpha))$ is equivalent to the event

$$
\left(\min_{i \in I_2} \mathcal{Z}^{i}_{\alpha} < \min_{i \in I_1 \cup I_3} \mathcal{Z}^{i} \right)
$$

\noindent and hence, by denoting

$$
qrum=\# \{\alpha, BU(\alpha) holds \}.
$$

\noindent Similarly to the previous steps, we may denote

\begin{equation} \label{ensembleEstimRum}
t(B)(\ell, \ell^{\prime})=\{ \sigma \in \mathcal{S}_n, \min_{i \in I_2} \sigma(i) < \min_{i \in I_1 \cup I_3} \sigma(i) \}.
\end{equation}

\noindent The probability law $simrum$ is still given by

\begin{equation*} \label{probaEstimRum}
\mathbb{P}(simrum =s)=\sum_{C\in B_{t}}\prod_{\alpha \in C}\mathbb{P}_{Z_{\alpha}}(t(B)(\ell, \ell^{\prime})),
\end{equation*}

\noindent for $ks \in \mathbb{N}$.\\

\noindent Before we conclude, let us address two points.\\

\noindent Point (a) The previous developments are based on choosing random permutations. This is very time-consuming when $n$ is large. Consider functions of the form 
$Z_{\alpha}(i)=iX_{\alpha}+Y_{\alpha}$ mod $n$, $1\leq i\leq n$, $1\leq \alpha \leq k$, is a way to quickly have almost permutations of a small number of repetitions among the set 
$\{Z_{\alpha}(i), \ 1\leq i \leq n\}$.\\

\noindent Point (b) To achieve the target property in Point (a), we may simply consider a random variable 
$$
Z_{\alpha}=(Z_{\alpha}(i), \ 1\leq i \leq n)
$$

\noindent with values in some space $D$, with a size at least equal to $n$. For $1 \leq i\neq j\leq n$, $1\leq \alpha \leq k$, denote the probability of the event that the two lines 
$\mathcal{Z}^{i}_{\alpha}$ and $\mathcal{Z}^{j}_{\alpha})$ have at least on common coordinate by

\begin{eqnarray}
\hspace{4cm}p_{i,j}&=&\mathbb{P}(\mathcal{Z}^{i}_{\alpha}=\mathcal{Z}^{j}_{\alpha}))  \notag\\
&=&\mathbb{P}_{Z_{\alpha}}(\{ \sigma \in \mathcal{S}_n, \sigma(i)=\sigma(j) \}) \notag
\end{eqnarray}

\noindent By independence and stationary, the probability of the complement of the event $D_{i,j}$ that the lines $\mathcal{Z}^{i}$ and $\mathcal{Z}^{j}$, $1\leq i \leq j \leq$, are disjoint is

$$
\mathbb{P}(D_{i,j})=\left(1 - p_{i,j}\right)^k.
$$

\noindent Next, denote by $\mathcal{C}_{n,r}$ the class of lines $1\leq i\neq j\leq n$ such that the lines $\mathcal{Z}^{i}$ and $\mathcal{Z}^{j}$ have exactly $r\geq 1$ common coordinates. Denote 
the probability of the event $D_{n}$ that all the lines are mutually disjoint each other. We have

$$
\mathbb{P}(D_n) = 1 - \left( \sum_{r=1}^{k} \#(\mathcal{C}_{n,r}) p_{i,j}^{r}(1-p_{i,j})^{k-r} \right),
$$

\noindent or

\begin{equation} \label{conditionEstimRum}
\mathbb{P}(D_n) = 1 - \left( \sum_{r=1}^{k} \#(\mathcal{C}_{n,r}) \mathbb{P}_{Z_{\alpha}}(\mathcal{D}_{i,j})^{r}(1-\mathbb{P}_{Z_{\alpha}}(\mathcal{D}_{i,j}))^{k-r} \right),
\end{equation}

\noindent with

$$
\mathcal{D}_{i,j} = \{ \sigma \in \mathcal{S}_n, \sigma(i)=\sigma(j) \}.
$$

\noindent \textbf{Conclusion}. We are now able to have a partial conclusion.\\ 

\noindent (a) If the probability in Formula (\ref{conditionEstimRum}), which is $\mathbb{P}(D_n)$ is zero, then set on which we compute the estimated similarity is given by
Formula (\ref{ensembleEstimRum}) and the probability law of the estimated similarity is given by Formula (\ref{probaEstimRum}).\\

\noindent (b) If the probability in Formula (\ref{conditionEstimRum}), which is $\mathbb{P}(D_n)$, is small enough, we approximate the probability law of the estimated similarity is given by Formula (\ref{probaEstimRum}).\\

\noindent (c) If we go back to the \textit{minhashing} functions and consider the functions $Z_{\alpha}$ as random permutations picked on the uniform sampling, the estimated similarity, on the base of Point (a) and (b), the exact mathematical expectation of any event $U_{\ell}(\alpha)=U_{\ell^{\prime}}(\alpha)$ for any $1\leq \alpha \leq k$. Why?\\  

\noindent In this pure and uniform random scheme, all the lines of $\mathcal{Z}$ are disjoint and the realization of the event

$$
t(B)(\ell, \ell^{\prime})=\{ \sigma \in \mathbb{S}_n, \min_{i \in I_2} \sigma(i) < \min_{i \in I_1 \cup I_3} \sigma(i) \},
$$

\noindent is a pure matter of combinatorics. To realize this event, we have to choose of the $\{\sigma(i), i\in I_2 \}$ to be the unity and we have 
$\#(I_2)=\#(S_{\ell} \cap S_{\ell^{\prime}})$ ways to do it. So the probability of having this is

$$
\frac{\#(S_{\ell} \cap S_{\ell^{\prime}})}{n}.
$$

\noindent At the arrival, for all $\alpha$, $1\leq \alpha\leq k$, the binary random variable $S(\alpha)$ that is equal one if $U_{\ell}(\alpha)=U_{\ell^{\prime}}(\alpha)$ and zero otherwise, is a Bernoulli random variable with parameter $sim=\#(I_2)/n$, the number we can uniformly choose a permutation $\sigma$ such that $1 \in \{\sigma(i), \ i\in I_2\}$. We have $k$ independent Bernoulli random variables. Besides, the random  variable

$$
S_{k}=\sum_{1\leq \alpha\leq k} S(\alpha)
$$

\noindent is $k$ times the estimated similarity \textit{simrum}. We may then apply the probability formulas :\\

\noindent (1) Moments :
$$
\mathbb{E}(simrum)=sim \ \ and  \ \ \mathbb{V}ar(simrum)=sim(1-sim)/k.
$$

\noindent (2) Tchebychev Inequality.
$$
\mathbb{P}(|simrum \- \ sim| >\lambda) \leq \frac{sim(1-sim)}{k\lambda^2}, \lambda>0.
$$

\noindent (3) Gaussian Approximation. If the similarity is non-zero, we have
$$
\sqrt{\frac{k}{sim(1-sim)}} \biggr( simrtm - sim\biggr) \rightarrow \mathcal{N}(0,1);
$$

\noindent which gives the approximated confidence interval of 95\% percent

$$
sim \in [simrum-1.96 \frac{\sqrt{simrum(1-simrum)}}{k}, \ simrum+1.96 \frac{\sqrt{simrum(1-simrum)}}{k}]
$$

\noindent Of course, we may have given more properties of the Binomial law and apply them to the similarity.

\noindent We achieved the result we were targeting. Nevertheless, we want to give preliminary results for the general.\\

\bigskip \noindent \textbf{D - General Case}.\\

\noindent We suppose that the $Z_{\alpha}$'s are independent observations a the random variable defined on \{1, \ldots, 2\} represented by $(Z(1),\ldots, Z(n))$ taking with values set $\mathcal{Z}$ which is finite and let $z_0$ its minimum member. We may use again the reasoning above. Given the event $D_n$ - the lines $\mathcal{Z}^{i}$'s disjoint - the random variable $S(\alpha)$ is one if and only if 
$$
z_0 \in \{Z_{\alpha}(i), \ i \in I_2\}
$$  

\noindent Let us denote by $\mathbb{P}_{D_n}$ the conditional probability on $D_n$, that is $\mathbb{P}_{D_n}(B)=\mathbb{P}_{D_n}(B/D_n)$ for all $B\subset \{ 1,2,\ldots,n\}$,  and by 
$\mathbb{E}_{D_n}$ the mathematical expectation with respect to $\mathbb{P}_{D_n}$.\\

\noindent We then have that $S_{k}$ follows a binomial law with parameter $k$ and

$$
p=\mathbb{P}_{D_n}(z_0 \in \{Z_{\alpha}(i), \ i \in I_2\})
$$

\noindent In this case, $simrum \rightarrow p$ in probability as $k\rightarrow +\infty$, since

\begin{eqnarray*}
	\hspace{2cm}\mathbb{P}(|simrum -p|>\lambda)&=&\mathbb{P}_{D_n}(|simrum -p|>\lambda \ / \ D_n) \mathbb{P}(D_n)\\
	&=& \mathbb{P}(D_n) \frac{p(1-p)}{k\lambda^2}, \ \lambda.
\end{eqnarray*}

\noindent With such a lay out, we will be able to find out possible other choices of the $Z_{\alpha}$'s we have very quickly while $\mathbb{P}(D_n)$ close to one as near as possible. It is thought that the choice $Z_{\alpha}(i)=iX_{\alpha}+Y_{\alpha}$ mod $n$ is justified. We will come  back to this.\\

\bigskip \noindent \textbf{E - Interesting remark}.\\

\noindent At the light of what precedes, we see that the \textit{RU} algorithm might have been done with maximum of the lines at the place of the minimum. The probability parameter of the Bernoulli random variables $S(\alpha)$ would the number of ways to uniformly chose a permutation $\sigma$ such that $n \in \{\sigma(i), \ i\in I_2\}$. From then, we proceed as above.

\section{Conclusions and perspectives} \label{sec4}
At the end, we unveil the validity of the \textit{RU} algorithm, justified the convergence of estimated similarity index, and completely described its probability law of this index in the pure and uniform scheme. But, we applied the method by using it on the set that is easily formed. The building of the modified set allows to gain a great amount. A formula describing the similarity indices obtained from the direct algorithm and the modified one allows to work with the second and, at the end of procedure, to find the first. Beyond the uniform scheme, we extended the method to general probability law and provided the limit of the estimated similarity. From there, we provided a way to have a reasonable estimation based on random variables that may be formed in short times, to the contrary to random permutations. The impact of using minhashing function should be evaluated in the frame developed here.

\end{document}